\documentclass[12pt,a4paper]{amsart}
\usepackage[utf8]{inputenc}
\usepackage[T1]{fontenc}
\usepackage{amssymb}
\usepackage{mathrsfs} 
\usepackage{amsmath}
\usepackage{amsthm}
\usepackage{amscd}
\usepackage{paralist} 
\setdefaultenum{\upshape(i)}{}{}{}
\usepackage[pdftex,colorlinks,linkcolor=black,filecolor=black,citecolor=black,urlcolor=blue]{hyperref}  

\allowdisplaybreaks

\newtheoremstyle{mythmstyle}%
{1.5\baselineskip}
{\baselineskip}
{\itshape}
{}
{\bf}
{}
{0pt}
{} 

\newtheoremstyle{mydefstyle}%
{1.5\baselineskip}
{\baselineskip}
{}
{}
{\bf}
{}
{0pt}
{} 

\newtheoremstyle{mypreuvestyle}%
{\baselineskip}
{\baselineskip}
{}
{}
{\em}
{}
{0pt}
{} 

\newif\ifmynonumberenvi\mynonumberenvitrue
\theoremstyle{mythmstyle}
\newtheorem{proclaimmythm}[equation]{} 
\newtheorem*{proclaimmythm*}{}
\newenvironment{proclaim}[2][*]{\ifx*#1\mynonumberenvitrue\begin{proclaimmythm*}{\bf#2.} \ignorespaces\else\mynonumberenvifalse\begin{proclaimmythm}{.\kern0.5em\bf#2.}\label{#1} \ignorespaces\fi}{\ifmynonumberenvi\end{proclaimmythm*}\else\end{proclaimmythm}\fi}

\theoremstyle{mydefstyle}
\newtheorem{proclaimmydef}[equation]{} 
\newtheorem*{proclaimmydef*}{}
\newenvironment{definition}[2][*]{\ifx*#1\mynonumberenvitrue\begin{proclaimmydef*}{\bf#2.}\else\mynonumberenvifalse\begin{proclaimmydef}{.\kern0.5em\bf#2.}\label{#1} \ignorespaces\fi}{\ifmynonumberenvi\end{proclaimmydef*}\else\end{proclaimmydef}\fi}

\newenvironment{definition*}[1]{\begin{proclaimmydef*}{\bf#1.} \ignorespaces}{\end{proclaimmydef*}}

\def\QEDbox{\hbox{\lower2.3pt\vbox{\hrule\hbox
   {\vrule\kern1pt\vbox{\kern1.7pt\hbox{$\scriptstyle
   QED$}\kern.6pt}\kern1pt\vrule}\hrule}}}
\def\QED{\hskip0.01em plus 40pt\null{} \null\nobreak\hfill
   \kern3pt\QEDbox} 
\newcommand\QEDici{\\\noalign{\vskip-\baselineskip\smash{\hbox to\linewidth{\vrule width0pt \hfill\global\QEDdejaplacetrue\QEDbox}}\vskip-\baselineskip}}
\newif\ifQEDdejaplace\QEDdejaplacefalse

\theoremstyle{mypreuvestyle}
\newtheorem*{proclaimmypreuve}{}

\newif\ifmynonumberequation\mynonumberequationtrue

\makeatletter
\numberwithin{equation}{section}
\newenvironment{moneq}[1][*]{\ifx*#1\mynonumberequationtrue\begin{equation*}\else\mynonumberequationfalse\begin{equation}\label{#1}\fi}{\ifmynonumberequation\end{equation*}\@ignoretrue\else\end{equation}\@ignoretrue\fi\ignorespaces}
\makeatother

\def\dcap_#1{\mathchoice{%
          {\textstyle\bigcap\limits_{#1}}}%
          {\underset{#1}\cap}%
          {\underset{#1}\cap}%
          {\underset{#1}\cap}}
\def\dcup_#1{\mathchoice{%
          {\textstyle\bigcup\limits_{#1}}}%
          {\underset{#1}\cup}%
          {\underset{#1}\cup}%
          {\underset{#1}\cup}}
\def\ddcap_#1^#2{\mathchoice{%
          {\textstyle\bigcap\limits_{#1}^{#2}}}%
          {\overset{#2}{\underset{#1}\cap}}%
          {\overset{#2}{\underset{#1}\cap}}%
          {\overset{#2}{\underset{#1}\cap}}}
\def\ddcup_#1^#2{\mathchoice{%
          {\textstyle\bigcup\limits_{#1}^{#2}}}%
          {\overset{#2}{\underset{#1}\cup}}%
          {\overset{#2}{\underset{#1}\cup}}%
          {\overset{#2}{\underset{#1}\cup}}}

\newcommand\bigrestricted{{\kern1pt\vrule height3.3ex depth1.7ex width0.6pt\kern1pt}}
\newcommand\caprestricted{{\kern1pt\smash{\vrule height1.7ex depth0.7ex width0.6pt}\relax\vrule width0pt depth0pt\kern2pt}}
\newcommand\midrestricted{{\kern1pt\vrule height1.7ex depth0.9ex width0.6pt\kern1pt}}
\newcommand\bmidrestricted{{\kern1pt\vrule height2.7ex depth1.7ex width0.6pt\kern1pt}}
\newcommand\NN{\mathbf{N}}
\newcommand\restricted{{\kern1pt\vrule height1.3ex depth0.5ex width0.6pt\kern1pt}}

\def\oversetalign#1\to#2{\mathbin{\smash{\overset{{\text{\rlap{\hss#1}}}}{#2}}}}
\def\oversettext#1\to#2{\mathbin{\smash{\overset{{\text{{#1}}}}{#2}}}}

\newcommand\norm[2][1]{\ifcase#1
   \left\Vert#2\right\Vert 
   \or\Vert#2\Vert 
   \or\bigl\Vert#2\bigr\Vert 
   \or\Bigl\Vert#2\Bigr\Vert 
   \or\biggl\Vert#2\biggr\Vert 
   \or\Biggl\Vert#2\Biggr\Vert 
   \else\Vert#2\Vert\fi} 

\newcommand\hominprod[4][1]{\ifcase#1
   \left\langle\homind{#4}#2,#3\right\rangle 
   \or\langle\homind{#4}#2,#3\rangle 
   \or\bigl\langle\homind{#4}#2,#3\bigr\rangle 
   \or\Bigl\langle\homind{#4}#2,#3\Bigr\rangle 
   \or\biggl\langle\homind{#4}#2,#3\biggr\rangle 
   \or\Biggl\langle\homind{#4}#2,#3\Biggr\rangle 
   \else\langle\homind{#4}#2,#3\rangle\fi} 

\newcommand\homsuperinprod[4][1]{\ifcase#1
   \supinprsym{#4}\left(#2,#3\right) 
   \or\supinprsym#4(#2,#3) 
   \or\supinprsym#4\bigl(#2,#3\bigr) 
   \or\supinprsym#4\Bigl(#2,#3\Bigr) 
   \or\supinprsym#4\biggl(#2,#3\biggr) 
   \or\supinprsym#4\Biggl(#2,#3\Biggr) 
   \else\supinprsym#4(#2,#3)\fi} 

\newcommand\inprod[3][1]{\ifcase#1
   \left\langle#2,#3\right\rangle 
   \or\langle#2,#3\rangle 
   \or\bigl\langle#2,#3\bigr\rangle 
   \or\Bigl\langle#2,#3\Bigr\rangle 
   \or\biggl\langle#2,#3\biggr\rangle 
   \or\Biggl\langle#2,#3\Biggr\rangle 
   \else\langle#2,#3\rangle\fi} 
   
\newcommand\inprodsym{\inprod\cdot\cdot}

\newcommand\inprodd[3][1]{\ifcase#1
   \left\langle\mkern-3mu\left\langle#2,#3
        \right\rangle\mkern-3mu\right\rangle 
   \or\langle\mkern-3mu\langle#2,#3
        \rangle\mkern-3mu\rangle 
   \or\bigl\langle\mkern-4mu\bigl\langle#2,#3
        \bigr\rangle\mkern-4mu\bigr\rangle 
   \or\Bigl\langle\mkern-6mu\Bigl\langle#2,#3
        \Bigr\rangle\mkern-6mu\Bigr\rangle 
   \or\biggl\langle\mkern-8mu\biggl\langle#2,#3
        \biggr\rangle\mkern-8mu\biggr\rangle 
   \or\Biggl\langle\mkern-10mu\Biggl\langle#2,#3
        \Biggr\rangle\mkern-10mu\Biggr\rangle 
   \else\langle\mkern-3mu\langle#2,#3
        \rangle\mkern-3mu\rangle\fi} 

\newcommand\supinprsym{\mathscr{S}}

\newcommand\superinprod[3][1]{\ifcase#1
   \supinprsym\left(#2,#3\right) 
   \or\supinprsym(#2,#3) 
   \or\supinprsym\bigl(#2,#3\bigr) 
   \or\supinprsym\Bigl(#2,#3\Bigr) 
   \or\supinprsym\biggl(#2,#3\biggr) 
   \or\supinprsym\Biggl(#2,#3\Biggr) 
   \else\supinprsym(#2,#3)\fi} 

\newcommand\contrfoper{{\iota}}

\newcommand\contrf[3][1]{\ifcase#1
   \contrfoper\left(#2\right)#3 
   \or\contrfoper(#2)#3 
   \or\contrfoper\bigl(#2\bigr)#3 
   \or\contrfoper\Bigl(#2\Bigr)#3 
   \or\contrfoper\biggl(#2\biggr)#3 
   \or\contrfoper\Biggl(#2\Biggr)#3 
   \else\contrfoper(#2)#3\fi} 

\newcount\familycount\familycount=0
\newcommand\nextfamily{\vskip0.3\baselineskip\noindent\advance\familycount by1{\textbf{Family \the\familycount. }}}

\newcommand\ah{{\hat a}}
\newcommand\alphah{{\hat\alpha}}
\newcommand\betah{{\hat\beta}}
\newcommand\bh{{\hat b}}
\newcommand\CA{\mathcal{A}}
\newcommand\CC{\mathbf{C}}
\newcommand\ch{{\hat c}}
\newcommand\Dense{\mathcal{D}}
\newcommand\eexp{\mathrm{e}}
\newcommand\End{\mathrm{End}}
\newcommand\extder{\mathrm{d}}
\newcommand\fracp[2]{\frac{\partial#1}{\partial#2}}
\newcommand\gammah{\hat\gamma}
\newcommand\Hilbert{\mathcal{H}}
\newcommand\ie{i.e.}
\newcommand\labelIntroSHS{B}
\newcommand\labelIntroSUR{C}
\newcommand\labelSHCP{A}

\newcommand\Liealg[1]{\mathfrak{#1}}
\newcommand\mapob{\kern0.75em}
\newcommand\mo{^{-1}}
\newcommand\myquote[1]{``#1''}

\newcommand\parity[1]{\vert#1\vert}
\newcommand\refmetnaam[2]{#1\ref{#2}}
\newcommand\RR{\mathbf{R}}
\newcommand\scirc{\,{\raise 0.8pt\hbox{$\scriptstyle\circ$}}\,}
\newcommand\shifttag[1]{\kern#1&\kern-#1}

\begin{document}

\author{Gijs M. Tuynman}

\title{The Super Orbit Challenge}

\address{Laboratoire Paul Painlev{\'e}, U.M.R. CNRS 8524 et D{\'e}partement de Math{\'e}matiques, Facult{\'e} des Sciences et Technologies, Universit{\'e} de Lille, 59655 Villeneuve d'Ascq Cedex, France}
\email{FirstName[dot]LastName[at]univ-lille[dot]fr}

\begin{abstract}
When using the generally adopted definition of a super unitary representation, there are lots of super Lie groups for which the regular representation is not super unitary. 
I propose a new definition of a super unitary representation for which all regular representations are super unitary. 
I then choose a particular super Lie group (of Heisenberg type) for which I provide a list of super unitary representations in my new sense, obtained by a heuristic super orbit method. 
The super orbit challenge is to find a well defined {super orbit method} that will reproduces more or less my list of super unitary representations (or explains why they should not appear). 

\end{abstract}

\thanks{This work was supported in part by the Labex CEMPI  (ANR-11-LABX-0007-01).}

\maketitle

\section{Introduction}

In order to make this paper as easily accessible as possible, I will interpret a super Lie group as a super Harish-Chandra pair $(G_o,\Liealg g)$, even though I prefer to interpret them as a supermanifold $G$ with a compatible group structure (in the sense of $\CA$-manifolds \cite{Tu04}). 
In a super Harish-Chandra pair $(G_o,\Liealg g)$, $\Liealg g=\Liealg g_0 \oplus \Liealg g_1$ is a super Lie algebra (over $\RR$) and $G_o$ an ordinary Lie group acting on $\Liealg g$ such that:
\begin{enumerate}[{\labelSHCP}1.]
\item
the Lie algebra of $G_o$ is (isomorphic to) $\Liealg g _0$;

\item
the action of $G_o$ preserves each $\Liealg g_\alpha$ (the action is \myquote{even});

\item
the restriction of the $G_o$ action to $\Liealg g_0$ is (isomorphic to) the adjoint action of $G_o$ on it Lie algebra.

\end{enumerate}
The generally accepted definition of a super unitary representation of a super Lie group $(G_o,\Liealg g)$ is the one that can be found (among others) in \cite[Def. 2, \S2.3]{CCTV:2006} and \cite{AllHilLau:2013}. 
One defines a super Hilbert space $(\Hilbert,\inprodsym, \supinprsym)$,as a graded Hilbert space $\Hilbert = \Hilbert_0\oplus \Hilbert_1$ with scalar product $\inprodsym$ and super scalar product $\supinprsym$ (a graded symmetric non-degenerate sesquilinear form) satisfying the following conditions:
\begin{enumerate}[\labelIntroSHS1.]
\item\label{IntroSHS1label}
$\inprod{\Hilbert_0}{\Hilbert_1}=0$;

\item\label{IntroSHS2label}
for all homogeneous $x,y\in \Hilbert$ we have $\superinprod xy = i^{\parity x}\cdot \inprod xy$.

\end{enumerate}
With these ingredients a super unitary representation of $(G_o,\Liealg g)$ on the super Hilbert space $(\Hilbert,\inprodsym, \supinprsym)$ then is a couple $(\rho_o,\tau)$ in which $\rho_o$ is an ordinary unitary representation of $G_o$ on the Hilbert space $\Hilbert$ and $\tau:\Liealg g\to \End\bigl(C^\infty(\rho_o)\bigr)$ an even super Lie algebra representation of $\Liealg g$ on $C^\infty(\rho_o)$, the space of smooth vectors for $\rho_o$ defined by
$$
C^\infty(\rho_o) = \{\,\psi\in \Hilbert \mid g\mapsto \rho(g)\psi \text{ is a smooth map }G\to\Hilbert\,\}
\mapob, 
$$
satisfying the conditions:
\begin{enumerate}[\labelIntroSUR1.]
\item
for each $g\in G_o$ the map $\rho_o(g)$ preserves each $\Hilbert_\alpha$ (the representation is \myquote{even});

\item
for each $X\in \Liealg g_0$ (the Lie algebra of $G_o$!) the map $\tau(X)$ is the restriction of the infinitesimal generator of $\rho_o\bigl(\exp(tX)\bigr)$ to $C^\infty(\rho_o)$;

\item
for each $X\in \Liealg g_\alpha$ the map $\tau(X)$ is graded skew-symmetric with respect to $\supinprsym$;

\item
for all $g\in G_o$ and all $X\in \Liealg g_1$ we have
$$
\tau(g\cdot X) = \rho_o(g)\scirc \tau(X) \scirc \rho_o(g\mo)
\mapob,
$$
where on the left we denote by $g\cdot X$ the action of $G_o$ on $\Liealg g$.

\end{enumerate}
Unfortunately, already for the most simple super Lie group $\RR^{0\vert 1}$, the $0\vert1$-dimensional abelian super Lie group for which the super Harish-Chandra pair is $(\{e\}, \{0\}\oplus\RR)$, the (left-) regular representation is not super unitary in the above sense. 
The representation space is the space of (smooth) functions $C^\infty(\RR^{0\vert1})$ of a single odd variable, \ie, isomorphic to $\CC^2$ via $f(\xi) = a_0+a_1\xi$ and the infinitesimal action is given by the operator $\partial_\xi$, \ie, by the matrix $(\begin{smallmatrix} 0 & 1 \\ 0 & 0\end{smallmatrix})$. As $\CC^2=\CC\oplus\CC$ is $1\vert1$-dimensional, there is no possible choice for a super Hilbert space structure on $\CC^2$ for which the regular representation is super unitary. 

In \cite{Tuynman:2017}\footnote{As \cite{Tuynman:2017} was too long for most journals, a shortened version without the sections on Berezin-Fourier decomposition will appear as \cite{Tuynman:2018}.} I proposed to change the definition of a super Hilbert space to a triple $(\Hilbert, \inprodsym, \supinprsym)$ by changing the condition \refmetnaam{\labelIntroSHS}{IntroSHS2label} to
\begin{enumerate}[\labelIntroSHS'1.]
\setcounter{enumi}{1}
\item
$\supinprsym$ is continuous with respect to the topology of $\Hilbert$ defined by $\inprodsym$.

\end{enumerate}
But remember, $\supinprsym$ is just a non-degenerate graded symmetric sesquilinear form, not necessarily even nor homogeneous. 
And then I proposed to change the definition of a super unitary representation on a super Hilbert space $(\Hilbert, \inprodsym, \supinprsym)$ as a triple $(\rho_o,\Dense,\tau)$ in which $\rho_o$ is an ordinary unitary representation of $G_o$ on the Hilbert space $\Hilbert$ and $\tau:\Liealg g\to \End\bigl(\Dense\bigr)$ an even super Lie algebra representation of $\Liealg g$ on $\Dense \subset C^\infty(\rho_o) \subset \Hilbert$, a dense graded subspace of $\Hilbert$ contained in the set of smooth vectors of the unitary representation $\rho_o$, satisfying the conditions:
\begin{enumerate}[\labelIntroSUR'1.]
\item
for each $g\in G_o$ the map $\rho_o(g)$ preserves each $\Hilbert_\alpha$ (the representation is \myquote{even});

\item
for each $X\in \Liealg g_0$ (the Lie algebra of $G_o$!) the map $\tau(X)$ is the restriction of the infinitesimal generator of $\rho_o\bigl(\exp(tX)\bigr)$ to $\Dense$;

\item
for each $X\in \Liealg g_\alpha$ the map $\tau(X)$ is graded skew-symmetric with respect to $\supinprsym$;

\item
for all $g\in G_o$ and all $X\in \Liealg g_1$ we have
$$
\tau(g\cdot X) = \rho_o(g)\scirc \tau(X) \scirc \rho_o(g\mo)
\mapob;
$$

\item
$\Dense\subset \Hilbert$ is maximal with respect to the four conditions above.

\end{enumerate}
I then showed that the left-regular representation of any connected super Lie group is super unitary in this new sense.\footnote{A slightly less far going modification of the notion of a super Hilbert space and an associated notion of a super unitary representation is proposed in \cite
{deGoursacMichel:2015}.} 
In particular for the simplest example of the $0\vert1$-dimensional super Lie group cited above, it suffices to take an odd super scalar product $\supinprsym$ instead of an even one as imposed by the standard definition. 

Now I think that rendering all regular representations super unitary is sufficient reason to justify my change of the definition of a super unitary representation, but my initial motivation comes from a heuristic super version of the orbit method. 
In \cite{Tuynman:2010}\footnote{The timeline of the official publications is different from the production timeline as can be seen from the arXiv dates.} I introduced the notion of a mixed symplectic form and I showed that coadjoint orbits of a super Lie group carry in a natural way such a mixed symplectic form. In \cite{Tuynman:2009} (see also \cite{Tuynman:2010/2}) I then showed that representations associated to orbits with a non-homogeneous symplectic form appear in the (Fourier-Berezin) decomposition\footnote{The same Fourier-Berezin decomposition technique was used in \cite{AllHilWur:2016} to decompose the regular representation of the $0\vert1$-dimensional super Lie group described above, using a family of representations depending on an odd parameter.} of the regular representation of an explicit example of dimension $4\vert4$, justifying the introduction of non-homogeneous symplectic forms. 
Now there seems to be a certain reluctance to accept the notion of non-even symplectic forms (see for instance \cite{AllHilWur:2016}) and my \myquote{justifying} paper \cite{Tuynman:2009} has a serious drawback: half of the used procedure is heuristic and no (super) Hilbert spaces are mentioned. 

I still have no satisfactory way to produce (by means of super geometric quantization of super symplectic manifolds with a non-even symplectic form) structures that might lead to super Hilbert spaces; for super Lie groups, I only have a systematic way to produce representations (essentially on spaces of smooth functions) associated to coadjoint orbits and polarizations, and then I have to invent by hand the (super) Hilbert space structure adapted to such a representation and I have to adapt by hand the dependence on odd parameters linked to the specific orbit. 
But now that I have a convenient notion of a super unitary representation, I will give, for a particular (Heisenberg like) super Lie group of dimension $3\vert3$, a list of super unitary representations in my new sense. 
I am convinced this will be the complete list of all inequivalent irreducible super unitary representations of this group, but (of course) I have no proof and I might be wrong. 
And then the challenge is to find a systematic way to obtain them via a well-defined orbit method. 
The interested reader will find another example in \cite{Tuynman:2009} (for which I have the same conviction) to test any super orbit method, although in that example no mention is made of any kind of notion of super unitary representation.

\section{A super Lie group and a list of super unitary representations}

As a super Lie group of dimension $3\vert3$ our example $G$ \myquote{is} $\RR^{3\vert3}$ with three global even coordinates $a,b,c$ and three global odd coordinates $\alpha, \beta, \gamma$ and group law given by the multiplication
\begin{align*}
(a, b, \alpha, \beta,c,\gamma) \cdot (\ah, \bh, \alphah, \betah, \ch, \gammah)
&
=
\bigl(a+\ah, b+\bh, \alpha+\alphah, \beta+\betah, 
\\&
\kern4em
c+\ch + \tfrac12(a\bh -b\ah-\alpha\betah-\beta\alphah), 
\\&
\kern4em
\gamma + \gammah + \tfrac12(a\betah - \beta\ah+b\alphah - \alpha\bh)\bigr)
\mapob.
\end{align*}
As a super Harish-Chandra pair $(G_o,\Liealg g)$ it is given by the standard Heisenberg group $G_o = \RR^3$ of dimension $3$ with group law
\begin{moneq}
(a, b, c) \cdot (\ah, \bh, \ch)
=
\bigl(a+\ah, b+\bh, 
c+\ch + \tfrac12(a\bh -b\ah)
\bigr)
\mapob.
\end{moneq}
The super Lie algebra $\Liealg g = \Liealg g_0 \oplus \Liealg g_1$ of dimension $3\vert 3$ with three even basis vectors $e_0,e_1,e_2$ and three odd basis vectors $f_0,f_1,f_2$ is described by the commutators
\begin{moneq}
{}
[e_1,e_2]=e_0 = [f_1,f_2]
\qquad,\qquad
[e_1,f_2] = f_0 = [e_2,f_1]
\mapob,
\end{moneq}
all others either $0$ or determined by graded skew-symmetry. 
It is a central extension of the abelian super group of dimension $2\vert2$ by a $1\vert1$-dimensional center; at the algebra level the center is generated by the vectors $e_0,f_0$. 
And finally the (adjoint) action of $G_o$ on $\Liealg g$ is given by
\begin{moneq}
\begin{aligned}
(a,b,c)\cdot e_0 &= e_0
\quad,&
(a,b,c)\cdot e_1 &= e_1 - b\,e_0
\quad,&
(a,b,c)\cdot e_2 &= e_2 + a\,e_0
\\
(a,b,c)\cdot f_0 &= f_0
\quad,&
(a,b,c)\cdot f_1 &= f_1 + b\,f_0
\quad,&
(a,b,c)\cdot f_2 &= f_2 + a\,f_0
\mapob.
\end{aligned}
\end{moneq}

Once we have the description of our super Lie group, we can provide our list of seven families of super unitary representations. 
However, instead of providing the unitary representation $\rho_o$ of $G_o$ and the infinitesimal representation $\tau$, I will give the integrated version $\rho$, which is a bona fide representation of the full super group $G$. The unitary representation $\rho_o$ is directly obtained by putting $\alpha=\beta=\gamma=0$ in the expression for $\rho$, and $\tau$ is obtained by computing the derivatives of $\rho$ with respect to the six variables $a,b,c,\alpha, \beta, \gamma$ at the point $(a,b,c,\alpha,\beta,\gamma) = \mathbf0$. 
For the third family this will be done explicitly.

\nextfamily
We start with a family of $1$-dimensional representations depending on two real parameters $k,\ell$ and two odd parameters $\kappa, \lambda$. 
Our graded Hilbert space is given by $\Hilbert = \CC \oplus \{0\}$ with scalar product and super scalar product $\superinprod\chi\psi = \inprod{\chi}{\psi} = \overline\chi\cdot \psi$.
And then the representation $\rho$ is given by
\begin{moneq}
\rho(a,b,\alpha,\beta,c,\gamma)\psi
=
\eexp^{i(a k+b \ell + \alpha\kappa + \beta\lambda)}\,\psi
\mapob.
\end{moneq}

\nextfamily
For this family the Hilbert space is $\Hilbert = L^2(\RR^2)\oplus \{0\}$ with its standard scalar product and super scalar product given by
\begin{moneq}
\inprod\chi\psi = \superinprod\chi\psi
=
\int \overline{\chi(x,y)}\,\psi(x,y)\ \extder x \,\extder y
\mapob.
\end{moneq}
On this Hilbert space we define a one-parameter family of representations $\rho$ depending on a nonzero odd parameter $\kappa$ by
\begin{moneq}
\bigl(\rho(a,b,\alpha,\beta,c,\gamma)\psi\bigr)(x, y)
=
\psi(x +  b, y + a) \, \eexp^{i\alpha\kappa x} \, \eexp^{i\beta\kappa y}
\,\eexp^{i (\gamma + \frac12(\beta a + b \alpha))\kappa}
\mapob.
\end{moneq}

\nextfamily
Here the graded Hilbert space is $\Hilbert$ is given by $\Hilbert = L^2(\RR)\oplus L^2(\RR)$, which I interpret as functions of one even variable $x$ and one odd variable $\xi$ according to
\begin{moneq}
(\psi_0, \psi_1)\in L^2(\RR) \oplus L^2(\RR)
\qquad\cong\qquad
\psi(x,\xi)=\psi_0(x)+\xi\psi_1(x) 
\mapob.
\end{moneq}
The scalar product $\inprodsym$ and the (odd) super scalar product $\supinprsym$ are given by
\begin{align*}
\inprod\chi\psi 
&
= 
\int \overline{\chi_0(x)}\,{\psi_0(x)} + \overline{\chi_1(x)}\, {\psi_1(x)} \ \extder x
\\
\superinprod\chi\psi
&
=
\int \overline{\chi_0(x)}\,{\psi_1(x)} + \overline{\chi_1(x)}\, {\psi_0(x)} \ \extder x
\mapob.
\end{align*}
On this Hilbert space we define a $1$-parameter family of representations $\rho$ depending on a nonzero real parameter $k$ by
\begin{moneq}
\bigl(\rho(a,b,\alpha,\beta,c,\gamma)\psi\bigr)(x,\xi)
=
\psi(x+k a,\xi-k \alpha) \, \eexp^{ibx} \, \eexp^{-i\beta\xi}
\, \eexp^{i k(c + \frac12(ab-\alpha\beta))}
\mapob.
\end{moneq}
This means that the unitary representation $\rho_o$ is given by
\begin{align*}
\bigl(\rho_o(a,b,c)\psi\bigr)(x,\xi)
&
=
\psi(x+k a,\xi)\, \eexp^{ix b}
\,
\eexp^{i k(c + \frac12ab)}
\end{align*}
and the super Lie algebra representation is given by
\begin{moneq}
\begin{aligned}
\tau(e_0)\psi 
&
= i k \,\psi
&\quad
\tau(e_1)\psi
&
=
k\,\fracp\psi{x}
&\quad
\tau(e_2)\psi
&
=
ix\,\psi
\\
\tau(f_0)\psi
&
=
0
&
\tau(f_1)\psi
&
=
-k\,\fracp\psi{\xi} 
&\quad
\tau(f_2)\psi
&
=
-i\xi\,\psi 
\mapob.
\end{aligned}
\end{moneq}

\nextfamily
For this family the Hilbert space is the same as for the third family. 
On it we define a two-parameter family of representations $\rho$ depending on a real parameter $k$ and a nonzero odd parameter $\kappa$ by
\begin{moneq}
\bigl(\rho(a,b,\alpha,\beta,c,\gamma)\psi\bigr)(x, \xi)
=
\psi(x + a, \xi - \alpha) \, \eexp^{ib(k+\xi\kappa)} \, \eexp^{i\beta\kappa x}
\,\eexp^{i (\gamma + \frac12(\beta a - b \alpha))\kappa} 
\mapob.
\end{moneq}

\nextfamily
For this family the Hilbert space is again the same as for the third family. 
On it we define a two-parameter family of representations $\rho$ depending on a nonzero real parameter $k$ and a nonzero odd parameter $\kappa$ by
\begin{moneq}
\begin{aligned}
\bigl(\rho(a,b,\alpha,\beta,c,\gamma)\psi\bigr)(x,\xi)
&
=
\psi(x + a,\xi - \alpha) \, \eexp^{ib(xk+\xi\kappa)} \, \eexp^{i\beta(x\kappa - \xi k)}
\\&
\kern5em
\eexp^{i (\gamma + \frac12(\beta a - b \alpha))\kappa} \, 
\eexp^{ik (c + \frac12(a b +\beta \alpha))}
\mapob.
\end{aligned}
\end{moneq}

\nextfamily
Here the graded Hilbert space is $\Hilbert = \CC^2 \oplus \CC^2$, which I interpret as functions of two odd variables $\xi$ and $\eta$ according to
\begin{moneq}
\bigl((\psi_0,\psi_{12})\oplus (\psi_1,\psi_2)\bigr) \in \CC^2 \oplus \CC^2
\quad\cong\quad
\psi(\xi,\eta) = \psi_0 + \xi\,\psi_1 + \eta\,\psi_2 + \xi\eta\,\psi_{12}
\mapob.
\end{moneq}
The standard scalar product $\inprodsym$ and the super scalar product $\supinprsym$ are given by
\begin{align*}
\inprod\chi\psi
&
=
\overline{\chi_0}\,\psi_0 + \overline{\chi_{12}}\,\psi_{12} + \overline{\chi_1}\,\psi_1 + \overline{\chi_2}\,\psi_2
\\
\superinprod\chi\psi
&
=
\overline{\chi_0}\,\psi_{12} + \overline{\chi_{12}}\,\psi_{0} + \overline{\chi_1}\,\psi_2 - \overline{\chi_2}\,\psi_1
\mapob.
\end{align*}
On this Hilbert space we define a three-parameter family of representations $\rho$ depending on two real parameters $k,\ell$ and a nonzero odd parameter $\kappa$ by 
\begin{moneq}
\bigl(\rho(a,b,\alpha,\beta,c,\gamma)\psi\bigr)(\xi, \eta)
=
\psi(\xi - \beta,\eta - \alpha) \, \eexp^{ia(\xi \kappa + k)} \, \eexp^{ib(\eta \kappa + \ell)} 
\,\eexp^{i(\gamma - \frac12(\beta a + b \alpha))\kappa }
\mapob.
\end{moneq}

\nextfamily
For this family the Hilbert space is given by $\Hilbert = L^2(\RR)^2 \oplus L^2(\RR)^2$, which I interpret as functions of one even variable $x$ and two odd variables $\xi,\eta$ according to
\begin{multline*}
\bigl((\psi_0,\psi_{12})\oplus (\psi_1,\psi_2)\bigr) \in L^2(\RR)^2 \oplus L^2(\RR)^2
\\
\qquad\cong\qquad
\psi(x,\xi,\eta) = \psi_0(x) + \xi\,\psi_1(x) + \eta\,\psi_2(x) + \xi\eta\,\psi_{12}(x)
\mapob.
\end{multline*}
The scalar product $\inprodsym$ and the super scalar product $\supinprsym$ are given by
\begin{align*}
\inprod\chi\psi
&
=
\int
\overline{\chi_0(x)}\,\psi_0(x) + \overline{\chi_{12}(x)}\,\psi_{12}(x) + \overline{\chi_1}(x)\,\psi_1(x) + \overline{\chi_2(x)}\,\psi_2(x) \ \extder x
\\
\superinprod\chi\psi
&
=
\int
\overline{\chi_0(x)}\,\psi_{12}(x) + \overline{\chi_{12}(x)}\,\psi_{0}(x) + \overline{\chi_1(x)}\,\psi_2(x) - \overline{\chi_2(x)}\,\psi_1(x) \ \extder x
\mapob.
\end{align*}
On this Hilbert space we define a three parameter family of representations $\rho$ depending on two nonzero real parameters $k,p$ and a nonzero odd parameter $\kappa$ by
\begin{moneq}
\begin{aligned}
\bigl(\rho(a,b,\alpha,\beta,c,\gamma)\psi\bigr)(x,\xi,\eta)
&
=
\psi(x + a - p b, \xi- \alpha, \eta -  \beta)
\eexp^{ib( xk + \xi \kappa + p\eta \kappa)} \, 
\\&
\kern4em
\eexp^{i\beta (x \kappa - \xi k)} \,
\eexp^{i(\gamma -p b \beta
+ \frac12 (\beta a - b \alpha))\kappa}\, 
\\&
\kern6em
\eexp^{ ik(c + \frac12(a b +\beta \alpha - pb^2))}
\mapob.
\end{aligned}
\end{moneq}

\section{Concluding remarks}

$\bullet$ 
All (super) Hilbert spaces are interpreted as spaces of functions on super spaces of the form $\RR^{p\vert q}$, or more precisely as spaces whose elements we can interpret as (smooth) functions of $q$ odd variables with values in the space of square integrable functions of $p$ real variables, \ie,  $\Hilbert = C^\infty\bigl(\RR^{0\vert q}; L^2(\RR^p)\bigr)$. 
It then turns out that in all cases the super scalar product $\superinprod\chi\psi$ is realised as the (translation invariant) Berezin-Lebesgue integral $\int_{\RR^{p\vert q}} \overline{\chi(m)}\,\psi(m)\ \extder m$. 

$\bullet$
All (infinitesimal) representations $\tau$ act by differentiation or multiplication, and as such these operators act on the space of smooth functions $C^\infty(\RR^{p\vert q})$. 
In all cases the unmentioned dense subspace $\Dense$ then is given by
\begin{moneq}
\Dense = \{\,\psi\in C^\infty(\RR^{p\vert q}) \mid \forall k\in \NN\ \forall X_1, \dots, X_k\in \Liealg g : \tau(X_1) \scirc \cdots \scirc \tau(X_k)\psi \in \Hilbert\,\}
\ .
\end{moneq}

$\bullet$
In most of the families of representations we have required some of the parameters to be nonzero. 
Not because these representations do not exist when they take the value zero, but because in those cases the representation will certainly not be irreducible. 

$\bullet$
In my way of thinking, the first family is associated to (coadjoint) orbits of dimension $0\vert0$, the third family is associated to orbits of dimension $2\vert 2$ with an even symplectic form, the families $2$, $4$ and $6$ are associated to orbits of dimension $2\vert 2$ with an odd symplectic form, and the families $5$ and $7$ are associated to orbits of dimension $2\vert 2$ with a non-homogeneous symplectic form. 
The seventh family is atypical as it is obtained by a polarization that is not of \myquote{maximal} dimension. 
For $p=0$ we recover (apart from the term $\partial_\eta$ in $\tau(f_2)$) the fifth family with the additional variable $\eta$. 

$\bullet$
The attentive reader will have noticed that I have not been completely honest, as the families 4--7 do not fit my description of a super unitary representation. 
In particular $\rho_o$ is not an ordinary unitary representation of the ordinary Lie group $G_o$, due to the presence of the (supposedly nonzero) odd parameter $\kappa$. 
On the other hand, apart from the fact that some of the parameters are odd, all these representations definitely have a \myquote{unitary} look, especially when one realises that the super scalar product is defined by integration of the product $\overline\chi\,\psi$ with respect to a translation invariant \myquote{measure.} 
As moreover all these families are obtained in the same way, I am sorely tempted to want to enlarge the definition of a super unitary representation even more in order to include all these families (e.g., linking the scalar product $\inprodsym$ directly to the super scalar product $\supinprsym$ in the spirit of a Krein space as in \cite{deGoursacMichel:2015}, but by dropping the condition that $\rho_o$ should preserve $\inprodsym$). 
Unfortunately I have not (as yet) a satisfactory way to do so. 

$\bullet$ And last but not least: all feed back will be appreciated.

\providecommand{\bysame}{\leavevmode\hbox to3em{\hrulefill}\thinspace}
\providecommand{\MR}{\relax\ifhmode\unskip\space\fi MR }
\providecommand{\MRhref}[2]{%
  \href{http://www.ams.org/mathscinet-getitem?mr=#1}{#2}
}
\providecommand{\href}[2]{#2}

\end{document}